\documentclass[11pt]{amsart}
\usepackage{amssymb,pstricks,pstcol,pst-plot}

\definecolor{Red}{cmyk}{0,1,1,0}
\definecolor{OrangeRed}{cmyk}{0,0.6,1,0}            
\definecolor{DarkBlue}{cmyk}{1,1,0,0.20}
\definecolor{Black}{cmyk}{0,0,0,1}
\definecolor{Violet}{cmyk}{0.79,0.88,0,0}
\definecolor{Purple}{cmyk}{0.45,0.86,0,0}
\definecolor{myblue}{rgb}{0.66,0.78,1.00}

\parskip=\smallskipamount
\numberwithin{equation}{section}

\newtheorem{theorem}{Theorem}[section]
\newtheorem{lemma}[theorem]{Lemma}
\newtheorem{corollary}[theorem]{Corollary}
\newtheorem{proposition}[theorem]{Proposition}

\theoremstyle{definition}
\newtheorem{definition}[theorem]{Definition}

\newtheorem{problem}[theorem]{Problem}
\newtheorem{remark}[theorem]{Remark}

\newcommand{\C}{\mathbb{C}}

\newcommand{\N}{\mathbb{N}}

\newcommand{\Z}{\mathbb{Z}}

\newcommand{\R}{\mathbb{R}}

\newcommand{\cA}{\mathcal{A}}

\newcommand{\cC}{\mathcal{C}}

\newcommand{\cO}{\mathcal{O}}

\newcommand{\cU}{\mathcal{U}}

\def\e{\epsilon}

\def\dim{{\rm dim}\,}

\def\nbd{neighborhood}

\def\di{\partial}
\def\dibar{\bar\partial}
\def\bs{\backslash}
\def\disc{\triangle}

\def\wt{\widetilde}

\begin{document}
\title[Approximation of holomorphic mappings]
{Approximation of holomorphic mappings on strongly pseudoconvex domains}
\author{Barbara Drinovec-Drnov\v sek \& Franc Forstneri\v c}
\address{Institute of Mathematics, Physics and Mechanics, 
University of Ljubljana, Jadranska 19, 1000 Ljubljana, Slovenia}
\email{barbara.drinovec@fmf.uni-lj.si}
\email{franc.forstneric@fmf.uni-lj.si}
\thanks{Research supported by grants P1-0291 and J1-6173, Republic of Slovenia.}

%
%    General info
%

\date{July 6, 2006} 

\begin{abstract}
Let $D$ be a relatively compact strongly pseudoconvex domain in a Stein manifold $S$.
We prove that for every complex manifold $Y$ the set $\cA(D,Y)$ 
of all continuous maps $\bar D\to Y$ which are holomorphic in $D$ is a 
complex Banach manifold, and every $f\in \cA(D,Y)$ can be approximated 
uniformly on $\bar D$ by maps holomorphic in open neighborhoods of $\bar D$
in $S$. Analogous results are obtained for maps of class $\cA^r(D)$, $r\in \N$.
We also establish the Oka property for sections of continuous or smooth fiber bundles 
over $\bar D$ which are holomorphic over $D$ and whose fiber enjoys 
the Convex approximation property (Theorem~\ref{main}).
\end{abstract}
\maketitle

%
%
%  INTRODUCTION
%
%
\section{Introduction}
Let $D$ be a relatively compact strongly pseudoconvex domain in 
a Stein manifold $S$, and let $Y$ be a complex manifold.
For $r\in\Z_+$ we denote by $\cC^r(\bar D,Y)$ the set of all $\cC^r$ maps 
$\bar D\to Y$ and by $\cA^r(D,Y)$ the set of all $f\in \cC^r(\bar D,Y)$ which are 
holomorphic in $D$. Set $\cA(D,Y)=\cA^0(D,Y)$,
$\cC^r(\bar D,\C)=\cC^r(\bar D)$, $\cA^r(D,\C)=\cA^r(D)$. 

In this paper we use the technique of holomorphic sprays to investigate 
the structure of the space $\cA^r(D,Y)$ and to obtain approximation 
theorems for holomorphic maps.
A {\em holomorphic spray} is a family of maps in a 
given space of maps (such as $\cA(D,Y)$), depending holomorphically on 
a parameter $t$ in an open set $P\subset\C^N$;  the domination 
condition means submersivity with respect to the parameter 
$t\in P$ at $t=0$ (Def.\ \ref{Spray}).

\begin{theorem}
\label{Banach-manifold}
Assume that $D$ is a relatively compact, strongly pseudoconvex domain with
$\cC^\ell$ boundary $(\ell\ge 2)$ in a Stein manifold. For each $r\in \{0,1,\ldots,\ell\}$ 
and for every complex manifold $Y$ the space $\cA^r(D,Y)$ is a complex Banach manifold.
\end{theorem}

Theorem \ref{Banach-manifold} follows from Corollary \ref{deformations}
to the effect that the set of all maps in $\cA^r(D,Y)$ 
which are $\cC^0$-close to a given $f_0\in \cA^r(D,Y)$
is isomorphic to the set of all $\cC^0$-small sections of the 
complex vector bundle $f_0^*(TY)\to\bar D$ which are holomorphic in $D$
and of class $\cC^r(\bar D)$. This gives a local chart 
around $f_0 \in \cA^r(D,Y)$ with values in the Banach space 
$\Gamma_\cA^r(D,f_0^*(TY))$ of all sections 
$\bar D\to f_0^*(TY)$ of class $\cC^r(\bar D)$ which 
are holomorphic in $D$, and it is easily seen that these charts are 
holomorphically compatible. This Banach manifold structure  is natural 
in the sense that the evaluation map $D \times \cA^r(D,Y)\to Y$ is holomorphic. 

For the disc $\disc=\{z\in\C\colon |z|<1\}$, the Banach manifold 
$\cA(\disc,Y)$ is locally modeled on the complex Banach space 
$\cA(\disc)^n$ with $n=\dim Y$. In view of the Oka-Grauert principle
for vector bundles of class $\cA(D)$, due to Leiterer \cite{Leiterer1},
the same holds whenever every continuous complex vector bundle over $\bar D$
is topologically trivial. 

Lempert \cite{Lempert} obtained analogous results 
for the space of all $\cC^k$ maps  from a compact 
$\cC^k$ manifold $V$ into a complex manifold $M$, showing that
$\cC^k(V,M)$ is a Banach manifold if $k\in \Z_+$, resp.\ a
Fr\'echet manifold if $k=\infty$. He also described a complex manifold
structure on the space of real analytic maps by considering
Stein complexifications. Like in our case, 
maps $f\colon V\to M$ near a given map $f_0$ correspond to 
$\cC^0$-small sections of the vector bundle $f_0^*TM$. 
(In this direction see also \cite{IS,IS2,LZ1,LZ2};
for the corresponding theory of mappings between smooth real manifolds 
see Palais \cite{Palais1,Palais2}.)

We now turn to holomorphic approximation theorems.

\begin{theorem} 
\label{approximation-maps}
Assume that $D\Subset S$ is a strongly pseudoconvex domain with $\cC^\ell$ boundary
$(\ell \ge 2)$ in a Stein manifold $S$, $Y$ is a complex manifold and $r\in \{0,1,\ldots, \ell\}$. 
Any map in $\cA^r(D,Y)$ can be approximated in the $\cC^r(\bar D,Y)$ topology by maps 
holomorphic in open neighborhoods of $\bar D$.
\end{theorem}

When $S=\C^n$, $Y=\C$ and $r=0$, this classical result follows from the 
Henkin-Ram\'irez integral kernel representation of functions in $\cA(D)$ 
\cite[p.\ 87]{Henkin,Kerzman,Ramirez,HL1}. 
Another  approach is to cover $bD$ by finitely many open charts in which 
$bD$ is strongly convex, approximate the given function in $\cA(D)$ 
by a holomorphic function in each of the charts and 
patch these local approximations into a global one by solving 
a Cousin problem with bounds. This method also works for $r>0$ 
by using a solution operator to the $\dibar$-equation
for $(0,1)$-forms with $\cC^r$ estimates 
(Range and Siu \cite{Range-Siu}, Lieb and Range \cite{Lieb-Range},
Michel and Perotti~\cite{Michel-Perotti}).
For approximation on certain weakly pseudoconvex domains 
see Fornaess and Nagel \cite{FN}.

The kernel approach  essentially depends on the linear structure
of the target manifold and does not seem amenable to generalizations.
In this paper we adapt the second approach mentioned above to
an arbitrary target manifold $Y$ by using the method of sprays.
A crucial ingredient is a result
from \cite{BDF} concerning gluing of pairs
of sprays on special configurations of domains called 
{\em Cartan pairs}. Using this gluing technique we show that any 
$f_0\in\cA^r(D,Y)$ is the core map of a dominating spray 
(Corollary \ref{sprays-exist2}). This immediately implies  
Theorem \ref{Banach-manifold}, and Theorem \ref{approximation-maps}
is then obtained by the bumping method alluded to above. 
We prove the analogous approximation result for sections 
of holomorphic submersions (Theorem \ref{approximation-sections}).
When $r\ge 2$, a different proof can be 
obtained by using the fact that the 
graph $G_f=\{(z,f(z))\colon z\in \bar D\} \subset S\times Y$ of any map
$f\in\cA^2(D,Y)$ admits a basis of open Stein neighborhoods in $S\times Y$ 
\cite[Theorem 2.6]{BDF}. 

For approximation of holomorphic maps 
from planar Jordan domains to (almost) complex manifolds 
see D.\ Chakrabarti \cite[Theorem 1.1.4]{Chak}, \cite{Chak2}.

For maps from domains in open Riemann surfaces we obtain the
following stronger result which extends Theorem 5.1 in \cite{BDF}.

\begin{theorem}
\label{Riemann-surfaces}
Let $D$ be a relatively compact domain with $\cC^2$ boundary in an open Riemann
surface $S$ and let $Y$ be a (reduced, paracompact) complex space. Every continuous  
map $\bar D\to Y$ which is holomorphic in $D$ can be approximated
uniformly on $\bar D$ by maps which are holomorphic in open
neighborhoods of $\bar D$ in $S$. The analogous result holds
for maps of class $\cA^r(D,Y)$.
\end{theorem}

The domains of the approximating maps in Theorems \ref{approximation-maps}
and \ref{Riemann-surfaces} must in general shrink to $\bar D$. 
For target manifolds $Y$ satisfying the following {\em Convex approximation property}
we also obtain global approximation results.

%
%
%  CAP
%
\begin{definition}
\label{CAP}
{\rm \cite[Def.\ 1.1]{ANN}}
A complex manifold $Y$ enjoys the {\em Convex Approximation Property} 
{\rm (CAP)} if any holomorphic map from a neighborhood of a compact convex 
set $K\subset \C^n$ $(n\in \N)$ to $Y$ can be approximated uniformly on $K$ 
by entire maps $\C^n\to Y$. 
\end{definition}

Recall that a compact set $K$ in a complex manifold $S$ is 
$\cO(S)$-convex if for every point $p \in S\bs K$ there is a
holomorphic function $g\in\cO(S)$ satisfying $|g(p)|>\sup_{z\in K} |g(z)|$.

The following is a consequence of Theorem \ref{approximation-maps}
and the main result of \cite{ANN}; for a more general
result see Corollary \ref{improved-Oka} below.

\begin{corollary}
\label{cor1}
Assume that $S$ is a Stein manifold and $D\Subset S$ is a strongly pseudoconvex
domain with $\cC^\ell$  boundary $(\ell\ge 2)$ whose closure $\bar D$ is
$\cO(S)$-convex. Let $r\in\{0,1,\ldots,\ell\}$ and let 
$f\colon S \to Y$ be a $\cC^r$ map which is holomorphic in $D$. 
If $Y$ enjoys {\rm CAP} then $f$ can be approximated in the 
$\cC^r(\bar D,Y)$ topology by holomorphic maps $S\to Y$ which are 
homotopic to $f$. 
\end{corollary}

We shall now extend Corollary \ref{cor1} to sections of certain 
holomorphic fiber bundles over compact strongly pseudoconvex domains.

\begin{definition}
\label{Ar-bundle}
Assume that $D\Subset S$ is domain with $\cC^\ell$ boundary  $(\ell\ge 2)$ 
in a Stein manifold $S$, $Y$ is a complex manifold and 
$r\in \{0,1,\ldots,\ell \}$. An {\em $\cA_Y^r(D)$-bundle} is a fiber bundle 
$h\colon X\to \bar D$ with fiber $h^{-1}(z)\simeq Y$ $(z\in\bar D)$ which is smooth of 
class $\cC^r(\bar D)$  and is holomorphic over $D$.  
More precisely, every point $z_0\in \bar D$ admits a relatively open neighborhood 
$U$ in $\bar D$ and a $\cC^r$ fiber bundle isomorphism 
$\Phi\colon X|_U =h^{-1}(U) \to U\times Y$ 
which is holomorphic over $U\cap D$. 
\end{definition}

We denote by $\Gamma^r(\bar D,X)$ the set of all sections 
$f\colon \bar D\to X$ of class $\cC^r$, and by $\Gamma_\cA^r(D,X)$ the set 
of all $f\in \Gamma^r(\bar D,X)$ which are holomorphic over $D$.

%
%
%   SIMPLIED VERSION OF THE MAIN THEOREM
%
%
\begin{theorem} 
\label{main}
{\em (The Oka property for sections of $\cA^r(D)$-bundles)}
Assume that $D\Subset S$ is a strongly pseudoconvex domain
with $\cC^\ell$ boundary $(\ell\ge 2)$ in a Stein manifold $S$  
and $h\colon X\to \bar D$ is an $\cA_Y^r(D)$-bundle for some $r\in\{0,1,\ldots, \ell\}$. 
If the fiber $Y=h^{-1}(z)$ $(z\in \bar D)$ enjoys {\rm CAP} then 
\begin{itemize}
\item[(i)] every continuous section $f_0\in\Gamma(\bar D,X)=\Gamma^0(\bar D,X)$ 
is homotopic to a section $f_1\in \Gamma_\cA^r(\bar D,X)$, and  
\item[(ii)] every homotopy $\{f_t\}_{t\in[0,1]} \in \Gamma(\bar D,X)$ with 
$f_0,f_1\in \Gamma_\cA^r(D,X)$ can be deformed with fixed ends to a homotopy 
in $\Gamma_\cA^r(D,X)$.
\end{itemize}
\end{theorem}

Theorem \ref{main} is proved in \S 6; a more precise result,
with approximation on certain compact sets in $\bar D$, 
is Theorem \ref{main-bis}. 
Theorem \ref{main} also applies if the base $\bar D$ is a compact complex 
manifold with Stein interior and smooth strongly pseudoconvex boundary;
such manifold is diffeomorphic to a strongly pseudoconvex 
domain in a Stein manifold by a diffeomorphism which is holomorphic in 
the interior \cite{Heunemann3,Ohsawa} (some loss of smoothness may occur). 
In this connection we mention Catlin's  boundary version \cite{Cat}
of the Newlander-Nirenberg integrability theorem \cite{NN}.

For fiber bundles over {\em open} Stein manifolds 
(without boundary) the corresponding result was obtained in \cite{ANN}
(see also \cite{FP1,Gromov,L2}). 
The classical case when $X\to S$ is a principal $G$-bundle, 
with fiber a complex Lie group $G$, is due to H.\ Grauert \cite{G2,G3}.
His results imply that the holomorphic classification of principal $G$-bundles 
over Stein spaces agrees with the topological classification. 
(See also the exposition of H.\ Cartan \cite{Cartan} and the papers
\cite{FRa1,FRa2,Fsurvey,HL2}.) A.\ Sebbar \cite{Sebbar} proved 
the analogous result for principal $G$-bundles which are smooth on the closure 
of a pseudoconvex domain $D\subset \C^n$ and holomorphic over $D$.

The CAP property of the fiber $Y$
is both necessary and sufficient for validity of a stronger version of
Theorem \ref{main} (i) with approximation on compact holomorphically convex
subsets (see Theorem \ref{main-bis} (i)).
In the absence of topological obstructions, CAP of $Y$ also implies extendibility of 
holomorphic maps from closed complex subvarieties in Stein manifolds to $Y$
\cite{FOURIER,L3}.  Among the conditions implying CAP are (from strongest to weakest):
\begin{itemize}
\item[(a)]  homogeneity under a complex Lie group (Grauert \cite{G2}), 
\item[(b)] the existence of a dominating spray on $Y$ (Gromov \cite{Gromov}), and 
\item[(c)] the existence of a finite dominating family of sprays \cite{SUBELL}.
\end{itemize}
Further examples of manifolds enjoying CAP can be found in \cite{FLEX,ANN}.

The paper is organized as follows. 
In \S 2 we recall the notion of a (dominating) holomorphic spray 
and some of the related technical tools developed in \cite{BDF}. 
A key ingredient is Proposition \ref{gluing-sprays} concerning gluing 
of holomorphic sprays on Cartan pairs (Def.\ \ref{Cartan-pair}).
In sections 3 and 4 we obtain some technical results on the existence
and approximation of sprays. Theorem \ref{Banach-manifold}
is a  consequence of Corollary \ref{TNT} which provides 
linearization of a neighborhood of any section of class $\cA^r(D)$. 
In \S 5 we prove Theorems \ref{approximation-maps} and \ref{Riemann-surfaces}.
In \S 6 we prove Theorem \ref{main-bis} which includes Theorem \ref{main}.

%
%
%  Sprays, gluing lemma
%
%
\section{Gluing sprays on Cartan pairs}
In this section $D$ is a relatively compact,
strongly pseudoconvex domain with $\cC^\ell$ boundary $(\ell\ge 2)$ in a Stein manifold $S$
and $h\colon X\to \bar D$ is either an $\cA^r_Y(D)$-bundle for some $r\in\{0,\ldots,\ell\}$,
or the restriction to $\bar D$ of a holomorphic submersion $\wt h\colon \wt X\to S$. 
(Our results also apply if $h\colon X\to \bar D$ is a 
smooth submersion which is holomorphic over $D$.)

Let $VT_x X =\ker dh_x$ denote the {\em vertical tangent space} at a point $x\in X$.

%
%
%  SPRAYS
%
%
\begin{definition}
\label{Spray}
Assume that $P$ is an open set in $\C^N$ containing the origin.
An {\em $h$-spray of class $\cA^r(D)$ with the parameter set $P$} is a map 
$f\colon \bar D\times P\to X$ of class $\cC^r$ which is holomorphic 
on $D\times P$ and satisfies $h(f(z,t))=z$ for all $z \in \bar D$ 
and $t\in P$. We call $f_0=f(\cdotp,0)$ the {\em central} (or {\em core}) section 
of the spray. A spray $f$ is  {\em dominating on a subset} 
$K \subset \bar D$ if the partial differential
\[
	\di_t|_{t=0} f(z,t) \colon T_0 \C^N=\C^N  \to VT_{f(z,0)} X
\]
is surjective for all $z\in K$; if this holds with $K=\bar D$ 
then $f$ is {\em dominating}. 

A (dominating) {\em 1-parametric $h$-spray of class $\cA^r(D)$ 
and parameter set $P$} is a continuous map $f\colon  [0,1]\times\bar D\times P \to X$
such that for each $s\in[0,1]$, $f^s=f(s,\cdotp,\cdotp)\colon \bar D\times P\to X$ 
is a (dominating) $h$-spray of class $\cA^r(D)$, and the derivatives 
of $f^s$ of order $\le r$ are continuous with respect to all variables
(including $s$).
\end{definition}

For a product fibration $h\colon X= \bar D\times Y \to \bar D$, $h(z,y)=z$,
we can identify an $h$-spray $\bar D\times P \to X$ 
with a {\em spray of maps} $f\colon \bar D\times P \to Y$ by 
composing with the projection $\bar D\times Y\to Y$, $(z,y)\to y$.
Such a spray is dominating on $K\subset\bar D$ if 
$\di_t|_{t=0} f(z,t) \colon T_0 \C^N \to T_{f(z,0)} Y$ is surjective
for all $z\in K$.

In applications the parameter set $P$ will be allowed
to shrink around $0\in\C^N$. If $f$ is dominating on a compact subset $K$ of 
$\bar D$ then by shrinking $P$ we may assume that
$\di_t f(z,t) \colon T_t \C^N \to VT_{f(z,t)} X$
is surjective for all $(z,t)\in K\times P$.

\begin{remark}
Our  definition of a spray is similar to Def.\ 4.1 in \cite{BDF},
but is adapted for use in this paper. Sprays were introduced to the
Oka-Grauert theory by Gromov \cite{Gromov} to obtain approximation and 
gluing theorems for holomorphic sections. In Gromov's terminology 
(which was also used in \cite{FP1,FP2,FP3}) a dominating spray 
on a complex manifold $Y$ is a holomorphic map $s\colon E\to Y$ 
from the total space of a holomorphic vector bundle
$\pi\colon E\to Y$ such that $s(0_y)=y$
and $ds_{0_y}(E_y)=T_y Y$ for every point $y\in Y$.
When $E=Y\times \C^N$, we get for any holomorphic map $f_0\colon D\to Y$
a dominating spray $f\colon D\times \C^N\to Y$ (in the sense
of Def.\ \ref{Spray}) by setting $f(z,t)=s(f_0(z),t)$. 
Global dominating sprays (with the parameter set $P=\C^N$) 
exist only rarely. Sprays whose parameter set is a 
(small) open subset of $\C^N$ have been called 
{\em local sprays}  in the literature on the Oka-Grauert problem.
\end{remark}

%
%
%  Definition of a Cartan pair
%
%
\begin{definition}
\label{Cartan-pair}
A pair of open subsets $D_0,D_1 \Subset S$ in a Stein manifold $S$
is said to be a {\em Cartan pair} of class $\cC^\ell$ $(\ell\ge 2)$ if 
\begin{itemize}
\item[(i)]  $D_0$, $D_1$, $D=D_0\cup D_1$ and $D_{0,1}=D_0\cap D_1$ 
are strongly pseudoconvex with $\cC^\ell$ boundaries, and 
\item[(ii)] $\overline {D_0\backslash D_1} \cap \overline {D_1\backslash D_0}=\emptyset$ 
(the separation property). 
\end{itemize}
We say that {\em $D_1$ is a convex bump on $D_0$} if, in addition to the above,
there is a biholomorphic map from an open neighborhood of $\bar D_1$ in $S$ onto 
an open subset of $\C^n$ $(n=\dim S)$ which maps $D_1$ and $D_{0,1}$
onto strongly convex domains in $\C^n$.
\end{definition}

The following proposition and its proof are crucial for all that follows.

%
%
%  GLUING SPRAYS
%
%
\begin{proposition}
\label{gluing-sprays}
{\em (Gluing of sprays)}
Let $(D_0,D_1)$ be a Cartan pair of class $\cC^\ell$ $(\ell\ge 2)$ 
in a Stein manifold $S$ and set $D=D_0\cup D_1$, $D_{0,1}=D_0\cap D_1$. 
Assume that $h\colon X\to \bar D$ is either an 
$\cA^r_Y(D)$-bundle for some $r\in\{0,\ldots,\ell\}$,
or the restriction to $\bar D$ of a holomorphic submersion  $\wt X\to S$.  

Given an $h$-spray $f\colon \bar D_0\times P_0\to X$ $(P_0\subset\C^N)$ 
of class $\cA^r(D_0)$ which is dominating on $\bar D_{0,1}$,
there is an open set $P\subset\C^N$, with $0\in P \Subset P_0$, satisfying the following. 
For every $h$-spray $f' \colon \bar D_1 \times P_0 \to X$ of class $\cA^r(D_1)$ 
which is sufficiently $\cC^r$ close to $f$ on $\bar D_{0,1} \times P_0$ 
there exists an $h$-spray $g\colon \bar D \times P \to X$ 
of class $\cA^r(D)$, close to $f$ in the $\cC^r$ topology 
on $\bar D_0\times P$, such that
$g_0=g(\cdotp,0)$ is homotopic to $f_0$ on $\bar D_0$
and is homotopic to $f'_0$ on $\bar D_1$.
If $f$ and $f'$ agree to order $m\in\N$ along $\bar D_{0,1}\times\{0\}$ 
then $g$ can be chosen to agree to order $m$ 
with $f$ along $\bar D_{0}\times\{0\}$ and with 
$f'$ along $\bar D_{1}\times\{0\}$.

If in addition $\sigma\subset \bar D_0$ is a common zero set of 
finitely many $\cA^r(D_0)$ functions and $\sigma\cap\bar D_{0,1}=\emptyset$ 
then $g$ can be chosen such that $g_0$ agrees with $f_0$ to a 
given finite order on~$\sigma$.

The analogous result holds for  1-parametric sprays.
\end{proposition}

Proposition \ref{gluing-sprays}  follows from
Proposition 4.3 in \cite{BDF}. For later reference we recall 
the main steps of the proof:

{\em Step 1:} 
Lemma 4.4 in \cite{BDF} gives a domain $P_1\Subset P_0$ containing 
the origin and a {\em transition map} $\gamma \colon\bar D_{0,1}\times P_1\to\C^N$ 
of class $\cA^r(D_{0,1}\times P_1)$, close to the map 
$\gamma_0(z,t)= t$ in the $\cC^r$ topology (the closeness depending on the
$\cC^r$ distance of $f'$ to $f$ on $\bar D_{0,1} \times P_0$) and satisfying  
\[
	f(z,t)=f'(z,\gamma(z,t)), \quad z\in\bar D_{0,1},\ t\in P_1.
\]

{\em Step 2:}
Let $P \Subset P_1$ be a domain containing the origin $0\in\C^N$.
If $\gamma$ is sufficiently $\cC^r$-close to $\gamma_0$ 
on $\bar D_{0,1}\times P_1$  then Theorem 3.2 in \cite{BDF} furnishes maps 
$\alpha \colon\bar D_0\times P\to\C^N$ and $\beta\colon\bar D_1\times P\to\C^N$, 
of class $\cA^r$ on their respective domains, satisfying 
\[
	\gamma(z,\alpha(z,t))=\beta(z,t), \quad z\in \bar D_{0,1},\ t\in P.
\]

{\em Step 3:} 
It follows that the sprays $f(z,\alpha(z,t))$ and $f'(z,\beta(z,t))$,
defined on $\bar D_0\times P$ resp.\ on $\bar D_1\times P$,
amalgamate into a spray $g\colon \bar D\times P\to Z$ with the stated
properties. (Compare with (4.4) in \cite{BDF}.)

%
%
%  Approximation theorems 
%
%
\section{Approximation of sprays on convex domains}
The following approximation result will be used in \S 4.

%
%
%  Approximation of a spray
%
%
\begin{lemma}
\label{approximation1}
Assume that $C\subset B$ is a pair of bounded strongly convex domains 
with $\cC^2$ boundaries in $\C^n$, $P\subset \C^N$ is an open set 
containing the origin, and $Y$ is a complex manifold.
Given a spray of maps $f\colon \bar C\times P\to Y$ of class $\cA^r(C,Y)$
$(r\ge 2)$ whose central map $f_0=f(\cdotp,0)$ extends to a map 
$\bar B\to Y$ of class $\cA^r(B,Y)$, there exist 
\begin{itemize}
\item[(i)]  open sets $P_0\supset P_1\supset P_2\supset\cdots$ 
in $\C^N$, with $P_0\subset P$ and $0\in P_j$ for all $j\in\Z_+$,
\item[(ii)] relatively open domains $\Omega_j\subset \bar B\times\C^N$ satisfying 
\[
	(\bar C\times \bar P_0)\cup (\bar B\times \bar P_j) \subset \Omega_j,
  \quad j=0,1,2,\ldots, 
\] 
\item[(iii)] a sequence of maps $F_j\colon \Omega_j\to Y$
of class $\cA^r(\Omega_j,Y)$ $(j\in\Z_+)$
\end{itemize}
such that $F_j(\cdotp,0)=f_0$ on $\bar B$ for all $j\in\Z_+$
and the sequence $F_j$ converges to $f$ in the $\cC^r$ 
topology on $\bar C\times \bar P_0$ as $j\to\infty$.

If in addition we are given a vector bundle map 
$L\colon \bar B\times\C^N \to TY$ of class $\cA^{r}(B)$
covering $f_0$ (i.e., such that the linear map 
$L_z\colon \{z\}\times \C^N\to T_{f(z,0)} Y$ 
is of class $\cC^r$ in $z\in \bar B$ and is holomorphic in $B$) and satisfying 
$L_z = \di_t f(z,t)|_{t=0}$ for $z\in\bar C$ then the sequence $F_j$ can be 
chosen such that, in addition to the above, 
$\di_t F_j(z,t)|_{t=0}=L_z$ for every $z\in\bar B$ and $j\in\Z_+$.

If $Y=\C^M$ then the above hold for all $r\in \Z_+$.
\end{lemma}

\begin{proof}
Consider first the case $Y=\C$. 
We may assume that $C$ contains the origin $0\in\C^n$ in its interior. 
Taylor expansion in the $t$ variable gives
\[
	f(z,t)=f_0(z)+ \sum_{j=1}^N g_j(z,t)\, t_j,\quad 
	(z,t)\in \bar C\times P 
\]
for some $g_{j}\in \cA^r(C\times P)$. 
For each $s<1$ the function $g^s_j(z,t)=g_j(sz,t)$ is holomorphic
in $\frac{1}{s} C\times P \supset \bar C\times P$; 
by choosing $s$ close to one we insure that the approximation 
is as close as desired in $\cC^r(\bar C\times P)$.
Fix an $s$, choose a polydisc $P_0 \Subset P$ containing $0\in\C^N$
and apply Runge's theorem to approximate $g^s_j$ 
in $\cC^r(\bar C\times \bar P_0)$ by an entire function $\wt g_j$. 
The function 
\[
	F(z,t)=f_0(z)+ \sum_{j=1}^N \wt g_j(z,t)\, t_j, 
	\quad (z,t) \in\bar B\times \C^N
\]	 
then approximates $f$ on $\bar C \times \bar P_0$ and it agrees
with $f_0$ on $\bar B\times \{0\}$. This gives a desired sequence
$F_j$ on the fixed domain $\Omega=\bar B\times P$. 

Assume in addition that $L\colon \bar B\times \C^N\to \C$ is a map 
of class $\cA^r(B\times \C^N)$ which is linear in the second variable and
satisfies $L_z= \frac{\partial}{\di t}|_{t=0} f(z,t)$ 
for $z\in\bar C$. By Taylor's formula we have  
\[
	f(z,t)=f_0(z)+ L_z (t) + \sum_{j,k=1}^N g_{j,k}(z,t)t_j t_k,\quad 
	(z,t)\in \bar C\times P,
\]
for some $g_{j,k}\in \cA^r(C\times P)$. Approximating each 
$g_{j,k}$ in $\cC^r(\bar C\times \bar P_0)$ by an entire function 
$\wt g_{j,k} \colon\C^n\times\C^N\to \C$ and setting
\[
	F(z,t)=f_0(z)+ L_z (t) + \sum_{j,k=1}^N \wt g_{j,k}(z,t) t_j t_k,
	\quad (z,t)\in \bar B\times \C^N, 
\]
gives the desired approximation.

The case $Y=\C^M$ follows by applying the above result to each component.

Assume now that $Y$ is an arbitrary complex manifold and that $r\ge 2$.
By \cite[Theorem 2.6]{BDF} the graph $\{(z,f_0(z))\colon z\in\bar B\}$ 
admits an open Stein neighborhood $U\subset \C^n\times Y$. 
Choose a proper holomorphic embedding
$\Phi\colon U\to\C^M$, an open neighborhood $V\subset \C^M$ of $\Phi(U)$
and a holomorphic retraction $\iota \colon V\to \Phi(U)$ onto $\Phi(U)$.
We apply the already proved approximation result to the map 
$(z,t)\to \wt f(z,t):= \Phi(z,f(z,t))\in\C^M$ to get a sequence
$\wt F_j$  satisfying the conclusion of Lemma \ref{approximation1}
with respect to $\wt f$. Let $pr_Y\colon \C^n\times Y\to Y$  denote 
the projection onto the second factor.
Assuming that $\wt F_j$  approximates $\wt f$ sufficiently closely
on $(\bar C\times \bar P_0) \cup (\bar B\times \{0\})$,
the latter set has an open neighborhood $\Omega_j$ on which the map 
$F_j=pr_Y\circ \Phi^{-1}\circ \iota\circ \wt F_j$ is defined,
and the resulting sequence $F_j$ satisfies the conclusion.
\end{proof}

\begin{remark}
We shall use Lemma \ref{approximation1} only with
$Y=\C^M$ and $r\in\Z_+$. The reason for restricting to the case $r\ge 2$ 
for a general $Y$ was that, at the time of this writing, it was not known 
whether the graph of any $\cA(B)$ map with values in $Y$ 
admits an open Stein neighborhood; this has been proved
in a subsequent publication of the second author 
\cite[Theorem 1.2]{FF:manifolds}, and therefore Lemma \ref{approximation1}
holds for all $r\in\Z_+$ and an arbitrary manifold $Y$.
\end{remark}

%
%
%  SECTION 4: Sprays around sections
%
%
\section{Linearization around a section}
We shall prove that each section of class $\cA^r(D)$ over a strongly 
pseudoconvex domain is the core of a dominating spray 
(Corollary \ref{sprays-exist2}). This is the key to all main results of the paper.
Applying Cartan's Theorem B for vector bundles of class $\cA^r(D)$
we then obtain an up to the boundary version of Grauert's tubular 
neighborhood theorem (Corollary \ref{TNT}), and 
Theorem \ref{Banach-manifold} easily follows.
The main step is provided by the following result obtained by 
combining Propositions \ref{gluing-sprays}  and \ref{approximation1}. 
(Compare with \cite[Lemma 4.2]{BDF}.)

\begin{proposition}
\label{sprays-exist}
Assume that $D$ is a relatively compact 
strongly pseudoconvex domain with $\cC^\ell$ boundary 
$(\ell\ge 2)$ in a Stein manifold $S$, $r\in\{0,1,\ldots,\ell\}$, and
$h\colon X\to \bar D$ is either an $\cA^r_Y(D)$-bundle 
(Def.\ \ref{Ar-bundle}) or the restriction to $\bar D$ of 
a holomorphic submersion $\wt h\colon \wt X\to S$. 
Given a section $f_0\in \Gamma^r_{\cA}(D,X)$ and 
a surjective complex vector bundle map 
$L\colon \bar D\times \C^N \to VT(X)|_{f_0(\bar D)}$ of class
$\cA^r(D)$ which covers $f_0$, there exist a domain 
$P\subset \C^N$ containing the origin and a (dominating)  
$h$-spray $f\colon \bar D\times P\to X$ of class $\cA^r(D)$ satisfying 
\begin{equation}
	f(z,0)=f_0(z),\quad \di_t|_{t=0} f(z,t) =L_z, \quad z\in \bar D. \label{L(z)}
\end{equation}

Furthermore, given a homotopy of sections $f^s_0\in \Gamma^r_{\cA}(D,X)$
$(s\in[0,1])$, covered by a homotopy of surjective complex vector bundle maps 
$L^s\colon \bar D \times \C^N \to VT(X)|_{f^s_0(\bar D)}$ 
which are holomorphic over $D$ such that
for $s=0,1$ the map $f^s_0$ is the central map of a spray $f^s$ 
over $\bar D\times P$ and $\di_t|_{t=0} f^s(z,t) =L^s_z$  $(z\in \bar D)$, 
there exist a domain $P_1\subset P$ containing the origin and a 
$1$-parametric spray $F\colon [0,1]\times \bar D\times P_1\to X$ of class $\cA^r(D)$
such that $F^s=F(s,\cdotp,\cdotp)$ agrees with $f^s$ on $\bar D\times P_1$ for $s=0,1$ and 
\begin{equation}
\label{L^s(z)}
	F^s(z,0)=f^s_0(z),\quad \di_t|_{t=0} F^s(z,t) =L^s_z, 
	\quad z\in \bar D,\ s\in[0,1]. 
\end{equation}
If $f^s_0$ and $L^s$ are smooth of order $k$ with respect to the
parameter $s\in [0,1]$ then $F^s$ may also be chosen smooth of order $k$
in $s$.
\end{proposition}

\begin{proof}
The conditions on $h\colon X\to \bar D$ imply that each point $x_0\in X$ 
admits an open neighborhood $W\subset X$ isomorphic to a product 
$U\times V$, where $U$ is a relatively open subset of $\bar D$
and $V$ is an open subset of a Euclidean space $\C^l$, 
such that in the coordinates $x=(z,w)$ $(z\in U,\ w\in \C^l$),
$h$ is the projection $(z,w)\to z$. Such coordinate neighborhoods 
in $X$ will be called {\em special}.  

By \cite[Lemma 12.3]{HL0} there exist strongly pseudoconvex 
domains $D_0\subset D_1\subset\cdots\subset D_m=D$  with $\cC^\ell$
boundaries such that $\bar D_0\subset D$, and
for every $j=0,1,\ldots,m-1$ we have $D_{j+1}=D_j\cup B_j$ 
where $B_j$ is a convex bump on $D_j$ (Def.\ \ref{Cartan-pair}).
Each of the sets $B_j$ may be chosen so small that $f_0(\bar B_j)$ is 
contained in a special coordinate neighborhood of $X$.
The essential ingredient in the proof is Narasim\-han's lemma on local convexification.

By Lemma 5.3 in \cite{FP1} there exists a dominating $h$-spray $f$ with core $f_0$ 
over a neighborhood of $\bar D_0$. We recall the main idea of the proof 
and show that $f$ can be chosen to satisfy $\di_t|_{t=0} f(z,t) =L_z$
for $z\in \bar D_0$. Let  $\{e_j\}_{j=1}^N$ be the standard basis of $\C^N$. Set 
\[
	L_j(f_0(z)):=L_z e_j \in VT_{f_0(z)}X, \quad j=1,\ldots, N. 
\]
Note that $f_0(D)$ is a closed complex submanifold of $X|_D=h^{-1}(D)$ 
and hence admits an open Stein neighborhood $\Omega\subset X|_D$ 
\cite{Siu}. Each $L_j$ is a holomorphic section
of the vertical tangent bundle $VT(X)$ on the submanifold $f_0(D)$ of $\Omega$, 
and by Cartan's Theorem B it extends to a holomorphic vertical
vector field on $\Omega$. Denote by $\theta^j_t$ its flow. The map
\[
	f(z,t_1,\ldots,t_N) =
	\theta_{t_N}^N\circ\cdots\circ\theta_{t_{2}}^2\circ \theta_{t_{1}}^1 (f_0(z)),
\]
which is well defined and holomorphic for $z$ in a neighborhood of $\bar D_0$ 
and for $t=(t_1,\ldots,t_N)$ in an open set $P\subset\C^N$ containing the origin,
is then a dominating spray satisfying $f(\cdotp,0)=f_0$ and
$\di_t|_{t=0} f(z,t) =L_z$.

To find a desired spray on $\bar D$ we perform a stepwise 
extension of $f$ over the convex bumps $B_0,\ldots,B_{m-1}$. At the $j$-th step 
we assume that we have a spray $\bar D_j\times P_j\to X$ with the required
properties, and we shall approximate it by a spray  $\bar D_{j+1}\times P_{j+1}\to X$
with a possibly smaller parameter set $0\in P_{j+1}\subset P_j$.
Since all steps are of the same kind, it suffices to explain the first step $j=0$. 

Applying Lemma \ref{approximation1} with the sets $C=D_0\cap B_0$, $B=B_0$ 
(which are strictly convex in local holomorphic coordinates) 
and $Y=\C^M$ (since $f_0(\bar B_0)$ is contained in a 
special coordinate chart of $X$), we find an open set 
$0\in P'\subset P_0$,  a relatively open set
$\Omega \subset \bar B\times \C^N$ containing 
$(\bar C\times \bar P_0)\cup (\bar B\times \bar P')$,
and a map $f' \colon \Omega \to \C^M$ of class $\cA^r(\Omega)$ 
which approximates $f$ in the $\cC^r$ topology on $\bar C\times P_0$, 
such that $f(z,0)=f'(z,0)$ and $\di_t|_{t=0}f(z,t)= \di_t|_{t=0}f'(z,t)$ 
for $z\in \bar C$. If the approximation is 
sufficiently close then Lemma 4.4 in \cite{BDF} 
furnishes a map $\gamma(z,t)=t+c(z,t)$ of class $\cA^r(C\times P')$  
in a smaller parameter set $0\in P'\Subset P_0$ 
which is $\cC^r$ close to $\gamma_0(z,t) = t$, 
$c(z,t)$ vanishes to second order at $t=0$ for all $z\in \bar C$, 
and $f(z,t)=f'(z,\gamma(z,t))$ for $(z,t)\in \bar C\times P'$. 

Applying \cite[Theorem 3.2]{BDF} to $\gamma$ on the Cartan pair $(D_0,B_0)$
we obtain a smaller parameter set $P_1\subset P'$ and maps 
$\alpha(z,t)=t+a(z,t)$ on $\bar D_0\times P_1$, $\beta(z,t)=t+b(z,t)$
on $\bar B_0 \times P_1$, of class $\cA^r$ and close to $\gamma_0(z,t)=t$ 
on their respective domains, such that $a(z,t)$ and $b(z,t)$ vanish
to second order at $t=0$ and 
$\gamma(z,\alpha(z,t))=\beta(z,t)$ for $z\in \bar C$ 
and $t\in P_1$. (See Step (ii) in the proof of 
Proposition \ref{gluing-sprays} above.)
Now $(z,t)\to f(z,\alpha(z,t))$ is a spray of class $\cA^r(D_0)$ 
which is $\cC^r$-close to $f$ and agrees with $f$ to second order 
at $t=0$, $f'(z,\beta(z,t))$ is a spray 
with the analogous properties on $\bar B_0\times P_1$, 
and by construction the two sprays agree on $\bar C \times P_1$; 
hence they define a spray satisfying (\ref{L(z)}) on the set
$\bar D_1=\bar D_0\cup \bar B_0$. 
After $m$ steps of this kind we obtain the first part of the Proposition.

We continue with the parametric case (this will only be used in the proof
of Theorem \ref{main}). Fix an $s\in [0,1]$. By the first part of the 
Proposition there exists an $h$-spray
$f^s\colon \bar D\times P\to X$ of class $\cA^r(D)$ satisfying 
\[
	f^s(z,0)=f^s_0(z),\quad \di_t|_{t=0} f^s(z,t)= L^s_z,\quad z\in \bar D.
\]
For $s=0,1$ we use the already given sprays. We wish to choose these sprays
to depend continuously (or smoothly) on the parameter $s$. To do this, we shall first 
use a fixed spray $f^s$ to find a solution in an open interval $I_s\subset \R$ 
around $s$, and finally we shall patch these solutions together.

Fix a number $u \in[0,1]$. Since $L^u$ is surjective, 
there is a direct sum splitting 
$\bar D\times\C^N=E\oplus G$, where $E$ and $G$ are vector bundles of 
class $\cA^r(D)$ and $E_z=\ker L^u_z$ for each $z\in \bar D$.
(We use Theorem B for $\cA(D)$-bundles, due to Leiterer \cite{Leiterer2},
and Heunemann's approximation theorem \cite{Heunemann1};
compare with the proof of Lemma 4.4 and the Appendix in \cite{BDF}.) 
We split the fiber variable on $\{z\}\times\C^N$ accordingly 
as $t=t'_z\oplus t''_z \in E_z\oplus G_z$. 
Note that $L^u\colon G\to VT(X)|_{f^u(\bar D)}$
is a complex vector bundle isomorphism of class $\cA^r(D)$.
By the inverse function theorem there is an open interval 
$I_u=(u-\delta,u+\delta) \subset\R$ such that for each 
$s \in I_u\cap [0,1]$ there exists a unique section 
$g_s\colon \bar D\to G$ of class $\cA^r(D)$ satisfying
\[
	f^u(z,0'_z\oplus g_s(z)) = f^s_0(z), \quad z\in\bar D. 
\]
It follows that the map
\[
	H^s(z,t)=f^u\bigl(z,t'_z\oplus (t''_z + g_s(z)) \bigr) 
\]
is a dominating 1-parametric spray with the core $f^s_0$
for $s\in I_u\cap [0,1]$. For $s=u$ we have 
$g_u=0$ and $H^u=f^u$.

It remains to adjust the $t$-differential of $H^s$ at $t=0$. 
Elementary linear algebra shows that, 
after shrinking the interval $I_u$ if necessary,
there exist for every $s\in I_u$ a unique complex 
vector bundle automorphism $A^s \colon G\to G$ 
and a unique complex vector bundle map $B^s\colon E\to G$,
both of class $\cA^r(D)$ and continuous in $s$, such that the map 
\[
	F^s(z,t) = H^s\bigl(z,t'_z  \oplus (B^s_z t'_z + A^s_z t''_z) \bigr),
	\quad s\in I_u \cap [0,1], 
\] 
is a 1-parametric spray satisfying 
$\di_t|_{t=0} F^s(z,t)= L^s_z(t)$ for all $z\in \bar D$
and  $s\in I_u \cap [0,1]$.

The above argument gives a finite covering of $[0,1]$ 
by intervals $I_j=[a_j,b_j]$ $(j=0,1,\ldots,m)$, where 
$a_0=0 < a_1 < b_0 < a_2 < b_1 < \cdots < b_m=1$,
and for each $j$ a 1-parametric spray 
$F_j \colon I_j\times\bar D\times P\to X$
satisfying the conclusion of the Proposition on $I_j$. 

It remains to patch the sprays $\{F_j\}_{j=0}^m$ into a 1-parametric
spray $F$ on $[0,1]$. Since each pair of adjacent intervals 
$I_j$, $I_{j+1}$ intersect in the segment $[a_{j+1},b_j]$
while each three intervals are disjoint, it suffices 
to explain how to patch $F_j$ and $F_{j+1}$ to a spray
over $I_j\cup I_{j+1}$. Choose a point $u\in (a_{j+1},b_j)$.
Using a decomposition $\bar D\times \C^N =E\oplus G$
as above, with $E=\ker \di_t|_{t=0} F_j(u,\cdotp,t)$,
the implicit function theorem gives a segment
$J_j =[\alpha_j,\beta_j] \subset (a_{j+1},b_j)$
with $\alpha_j<u<\beta_j$, a polydisc $P_0\subset P$ containing 
$0\in\C^N$, and for each $s\in J_j$ a unique map 
$\gamma(s,z,t)=t'_z\oplus (t''_z + c(s,z,t))$ 
of class $\cA^r(D\times P_0)$ whose derivatives of order $\le r$ 
in $(z,t)\in \bar D\times P_0$ are continuous in all variables 
and $c(s,z,t)$ vanishes to second order at $t=0$, such that
\[
	F_j(s,z,\gamma(s,z,t)) = F_{j+1}(s,z,t), \quad 
	s\in I_j,\ z\in\bar D,\ t\in P_0.  
\]
(The special form of $\gamma$ is insured by the 
fact that the first order jets of the sprays $F_j$ and $F_{j+1}$ 
with respect to $t$ agree at $t=0$ for every $s\in I_j$. 
For more details see Lemma 4.4 in \cite{BDF}.)
Choose a smooth function $\chi\colon \R\to[0,1]$
such that $\chi(s)=0$ for $s\le \alpha_j$
and $\chi(s)=1$ for $s\ge \beta_j$.
The map 
\[
	(s,z,t) \mapsto F_j\bigl(s,z, t'_z\oplus (t''_z + \chi(s) c(s,z,t))\bigr) 
\]
is then a 1-parametric spray of class $\cA^r(D)$ 
satisfying the Proposition on $I_j\cup I_{j+1}$.
After $m$ steps we obtain a solution on $[0,1]$.
\end{proof}

A map $L$ as in Proposition \ref{sprays-exist} always exists for a sufficiently 
large $N$ as follows from Cartan's Theorem A for coherent sheaves of $\cA^r$ modules 
on strongly pseudoconvex domains (\cite[Theorem 6]{Heunemann2}, \cite{Lieb-Range}).
Hence we get

\begin{corollary}
\label{sprays-exist2}
{\em (Existence of dominating sprays)}
Given $h\colon X\to \bar D$ and a section $f_0\in \Gamma^r_{\cA}(D,X)$
as in Proposition  \ref{sprays-exist}, there exist a domain 
$0\in P\subset \C^N$ for some $N\in \N$ and a dominating $h$-spray 
$f\colon \bar D\times P\to X$ of class $\cA^r(D)$ with the central section $f_0$.
Furthermore, for any homotopy of sections $f_0^s \in \Gamma^r_{\cA}(D,X)$
$(s\in[0,1])$ there exists a homotopy of dominating $h$-sprays
$f^s\colon \bar D\times P\to X$ of class $\cA^r(D)$ such that
the core of $f^s$ equals $f^s_0$ for every $s\in[0,1]$.
\end{corollary}

Corollary \ref{sprays-exist2}  implies the following `up to the boundary' 
version of Grau\-ert's tubular neighborhood theorem which allows linearization 
of analytic problems near a given section.

\begin{corollary}
\label{TNT}
{\em  (Linearizing a neighborhood)}
Let $h\colon  X\to\bar D$ be as in Proposition \ref{sprays-exist}.
Given a section $f\in \Gamma_\cA^r(D,X)$ 
there exist a holomorphic vector bundle $\pi\colon E\to \wt D$ 
over an open neighborhood $\wt D\subset S$ of $\bar D$,  
a relatively open neighborhood $\Omega$ of the zero section 
of the restricted bundle $E|_{\bar D}:=\pi^{-1}(\bar D)$ and a fiber 
preserving $\cC^r$ diffeomorphism $\Phi\colon \Omega\to\Phi(\Omega) \subset X$
which is biholomorphic on $\Omega\cap \pi^{-1}(D)$ and 
maps the zero section of $E|_{\bar D}$ onto $f(\bar D)$.
\end{corollary}

\begin{proof}
Let $f\colon \bar D\times P\to X$ be a dominating spray 
furnished by Corollary \ref{sprays-exist2}.
There is a splitting $\bar D\times \C^N = E\oplus E'$ of $\cA^r(D)$-vector bundles,
where $E'_z=\ker \di_t|_{t=0} f(z,t)$ and $E$ is a complementary bundle
(Theorem B for $\cA^r$ bundles, \cite{Heunemann2}).
The restriction of $f$ to $\Omega= E \cap P$ satisfies  
Corollary \ref{TNT}.
\end{proof}

Note that $E$ in Corollary \ref{TNT} is just the normal bundle of 
the given section $f\colon \bar D\to X$, and it can be identified with the 
restriction of the vertical tangent bundle $VT(X)=\ker dh$ to $f(\bar D)$. 
When $r\ge 2$ and $h$ extends to a holomorphic submersion $\wt X\to \wt D$ 
over an open neighborhood $\wt D$ of $\bar D$ in $S$, $f(\bar D)$ admits a basis of open Stein 
neighborhoods in $\wt X$ \cite[Theorem 2.6]{BDF} and the conclusion of 
Corollary \ref{TNT} easily follows from standard methods.

Corollary \ref{TNT} implies the following result concerning 
the deformation space of a map in $\cA^r(D,Y)$.
For related results concerning maps in $L^2$-Sobolev classes from 
certain bordered Riemann surfaces to (almost) complex manifolds
see Ivashkovich and Shevchishin \cite{IS,IS2}.

\begin{corollary}
\label{deformations}
{\em (The deformation space of an $\cA^r$ map)}
Assume that $D$ is a relatively compact domain with strongly pseudoconvex 
boundary of class $\cC^\ell$ $(\ell\ge 2)$ in a Stein manifold. 
If $r\in\{0,1,\ldots,\ell\}$ and $Y$ is an arbitrary complex manifold
then for any $f_0\in \cA^r(D,Y)$ the space of all  
$f\in\cA^r(D,Y)$ which are sufficiently $\cC^0$-close to $f_0$ is 
isomorphic to the space of $\cC^0$-small sections in $\Gamma_{\cA}^r(D,E)$, 
where $E=f_0^*(TY)$.
\end{corollary}

\begin{proof}
We may consider maps $f\colon \bar D\to Y$ as sections of the product fibration
$h\colon X=\bar D\times Y \to \bar D$. Fix $f_0\in \Gamma_\cA^r(D,X)$. 
Let $\Phi\colon\Omega\to\Phi(\Omega)\subset X$ be as in Corollary \ref{TNT}, 
where $\Omega$ is an open neighborhood of the zero section 
in the complex vector bundle $E=f_0^*(TY)$ and $\Phi$ maps the zero section
of $E$ onto $f_0(\bar D)$. If $f\colon \bar D\to X$ is a section in 
$\Gamma_\cA^r(D,X)$ which is sufficiently uniformly close to $f_0$
then $f(\bar D)\subset \Phi(\Omega)$, and hence $f=\Phi\circ f'$ 
for a unique  $f'\in \Gamma_{\cA}^r(D,E)$
with $f'(\bar D)\subset \Omega$.
\end{proof}

{\em Proof of Theorem \ref{Banach-manifold}.}
The proof is similar to Lempert's construction in \cite[\S 2]{Lempert},
and our Corollary \ref{deformations} plays a similar role 
as Lemma 2.1 in that paper. (For the smooth case see also Palais
\cite{Palais1,Palais2}.) The key difference is that the sprays 
in our paper must be holomorphic in the base variable
$z\in D$, and their construction is one of the main technical 
ingredients of our proof.

Given $f_0\in \cA^r(D,Y)$, Corollary \ref{deformations}
furnishes a contractible local chart $\cU\subset \cA^r(D,Y)$ 
consisting of all maps $z\to \Phi(z,\xi(z))$ $(z\in\bar D)$,
where $\xi \in \Gamma_{\cA}^r(D,f_0^*(TY))$
is a section with range in an open, fiberwise 
contractible set $\Omega_0 \subset E_0=f_0^*(TY)$
containing the zero section. 
The transition map between any pair of such charts is of the form 
\[
	\Gamma_{\cA}^r(D,E_0) \ni \xi \to  \wt\xi \in \Gamma_{\cA}^r(D,E_1),
\]
where $\wt \xi (z)= \Psi(z, \xi(z))$ $(z\in\bar D)$
for a fiber preserving diffeomorphism $\Psi$ of class $\cA^r$
from an open subset of $\Omega_0\subset E_0$ 
onto an open subset of $E_1=f_1^*(TY)$. 
(Actually $E_0$ and $E_1$ are isomorphic since $f_1$ is isotopic to $f_0$
due to fiber contractibility of the set $\Omega_0$.)
In local coordinates $(z,w)$ on $E_0$ resp.\ $E_1$, 
$\Psi$ is of the form $\Psi(z,w)=(z,\psi(z,w))$.
The differential of the transition map $\xi\to \wt \xi$ at $\xi_0$, 
applied to a tangent vector $\eta\in\Gamma_\cA^r(D,E_0)$, equals
\[
		z\to \di_2 \Psi(z,\xi(z))\, \eta(z), \quad z\in\bar D.
\] 
Since the partial differential $\di_2 \Psi$ is nondegenerate
and the $\cC^r$ regularity up to the boundary is preserved
when differentiating on the fiber variable,
we see that the differential is a complex Banach space isomorphism
effected by $\di_2 \Psi(\cdotp,\xi(\cdotp))$, 
and hence the transition map is biholomorphic. 
\qed

%
%
%  SECTION 5
%
%
\section{Uniform approximation of holomorphic sections}
In this section we prove the following approximation theorem
which includes Theorem \ref{approximation-maps} as 
a special case.

\begin{theorem}
\label{approximation-sections}
Assume that $D\Subset S$ is a strongly pseudoconvex domain  
of class $\cC^\ell$ $(\ell\ge 2)$ in a Stein manifold $S$.
Let $h\colon X\to S$ be a holomorphic submersion of a complex manifold $X$
onto $S$ and let $r\in\{0,1,\ldots,\ell\}$.
Every section $f\in \Gamma_\cA^r(D,X)$ can be approximated in the 
$\cC^r$ topology on $\bar D$ by sections which are holomorphic in 
open neighborhoods of $\bar D$ in $S$.
\end{theorem}

\begin{proof}
Let $n=\dim S$ and $n+m=\dim X$.
Denote by $pr_1$ resp.\ $pr_2$ the coordinate projections of $\C^n\times\C^m$
onto the respective factors $\C^n$, $\C^m$, and let $B\subset\C^n$, $B'\subset \C^m$
denote the unit balls.

Since $h\colon X\to S$ is a holomorphic submersion,  there exist
for each point $x_0\in X$ open neighborhoods $x_0\in W\subset X$,
$h(x_0)\in V\subset S$, and biholomorphic maps 
$\phi\colon V\to B\subset \C^n$,
$\Phi \colon W\to B\times B'\subset \C^n\times\C^m$,
such that $\phi(h(x))=pr_1(\Phi(x))$ for every $x\in W$.
Thus $\Phi$ is of the form
\[
	\Phi(x)= (\phi(h(x)),\phi'(x)) \in B\times B', \quad x\in W,
\]
where $\phi'=pr_2\circ\Phi$. Let us call such $(W,V,\Phi)$ 
a {\em special coordinate chart} on $X$. Note that $h(W)=V$.

Fix a section $f\colon \bar D\to X$ in $\Gamma_{\cA}^r(D,X)$.
Using Narasimhan's lemma on local convexification of a 
strongly pseudoconvex hypersurface 
we find finitely many special coordinate charts 
$(W_j,V_j,\Phi_j)$ on $X$, with $\Phi_j=(\phi_j\circ h,\phi'_j)$,
such that $bD\subset \cup_{j=1}^{j_0} V_j$ 
and the following hold for $j=1,\ldots,j_0$:
\begin{itemize}
\item[(i)]   $\phi_j(bD\cap V_j)$ is a strongly convex hypersurface 
in the ball $B\subset \C^n$, 
\item[(ii)]  $f(\bar D\cap V_j) \subset W_j$, and 
\item[(iii)] $\overline{\phi'_j(f(D\cap V_j))} \subset B'$.
\end{itemize}

Choose a number $c<1$ sufficiently close to $1$ such that the sets
$U_j= \phi_j^{-1}(cB) \Subset V_j$ $(j=1,\ldots,j_0)$ still cover $bD$.

By a finite induction we shall construct an increasing sequence of
strongly pseudoconvex domains with $\cC^\ell$ boundaries
$D=D_0\subset D_1\subset\cdots\subset D_{j_0} \Subset S$
and sections $f_k\in \Gamma^r_\cA(D_k,X)$ $(k=0,1,\ldots,j_0)$,
with $f_0=f$, such that for every $k=1,\ldots, j_0$ the restriction 
of $f_k$ to $\bar D_{k-1}$ will be close to the previous section $f_{k-1}$
in $\Gamma_{\cA}^r(\bar D_{k-1},X)$. The domain $D_k$ 
of $f_k$ will in general depend on the rate of approximation
on $D_{k-1}$ and will be chosen such that 
\[
	D_{k-1}\subset D_k \subset D_{k-1}\cup V_k, \quad bD_{k-1} \cap U_k\subset D_k
\]
for $k=1,\ldots,j_0$; that is, we enlarge $D_{k-1}$ 
inside the $k$-th coordinate neighborhood $V_k$ so that
the part of $bD_{k-1}$ which lies in the smaller set $U_k$ 
is contained in the next domain $D_k$. 
The final domain $D_{j_0}$ will contain $\bar D$ in its interior, 
and the section $f_{j_0} \in \Gamma_\cA^r(D_{j_0},X)$ will approximate 
$f$ as close as desired in $\cC^r(\bar D,X)$. 
To keep the induction going we will also insure at every step
that the properties (ii) and (iii) above remain valid with 
$(D,f)$ replaced by $(D_k,f_k)$ for all $j=1,\ldots,j_0$.

Since all steps will be of exactly the same kind, we shall explain
how to get $(D_1,f_1)$ from $(D,f)=(D_0,f_0)$. We begin by finding a
domain $D'_1$ with $\cC^\ell$ boundary in $S$ which is 
a convex bump on $D=D_0$ (Def.\ \ref{Cartan-pair})
such that $\overline  U_1 \cap \bar D_0 \subset D'_1 \Subset V_1$.
We shall first find a set $\wt D'_1 \Subset B$ with desired properties 
and then take $D'_1=\phi_1^{-1}(\wt D'_1)$.

Choose a smooth function $\chi\ge 0$ with compact support contained 
in $B\subset\C^n$  such that $\chi=1$ on $cB$. 
Let $\tau\colon B\to\R$ be a strongly convex defining
function for the domain $\phi_1(D \cap V_1) \subset B$.
Choose $c'\in (c,1)$ close to $1$ such that the hypersurface
$\phi_1(bD\cap V_1)=\{\tau=0\}$ intersects the sphere 
$\{|z|=c'\}$ transversely. If $\delta >0$ is chosen sufficiently small 
then the set  
\[
	\{z\in \C^n \colon |z| < c',\ \tau(z) < \delta \chi(z) \}
\]
satisfies the required properties, except that it is not
smooth along the intersection of the hypersurfaces $\{|z| = c'\}$
and $\{\tau = \delta \chi\}$. By rounding off the corners
of the intersection we get a strongly convex set $\wt D'_1$
in $B$ such that $D'_1=\phi_1^{-1}(\wt D'_1)\subset V_1$ satisfies the
desired properties (Fig.\ \ref{Fig1}).  

%
%
%  Figure 1: The domains $D'_1$ and $D_1$
%
%
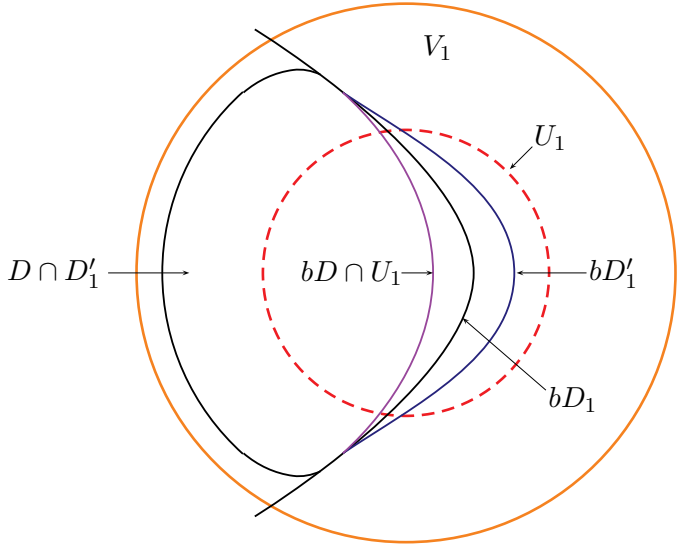
\begin{figure}[ht]
\psset{unit=0.6cm, linewidth=0.7pt}  
\begin{pspicture}(-8,-6)(8,6)

\pscircle[linecolor=OrangeRed,linewidth=1pt,fillstyle=none](0,0){6}
\pscircle[linecolor=Red,linewidth=1pt,linestyle=dashed](0,0){3.2}

\pscurve(-3.35,5.4)(1.5,0)(-3.35,-5.4)                               % The boundary of $D_1$
\psecurve[linecolor=DarkBlue](-4,5)(-1.4,4)(2.4,0)(-1.4,-4)(-4,-5)   % The boundary of $D'_1$
\psecurve[linecolor=Purple](-4,5)(-1.4,4)(0.6,0)(-1.4,-4)(-4,-5)     % The boundary of $D$
\psarc(0,0){5.4}{132}{228}                                           % arc - boundary of $D'_1\cap D$                         
\psecurve(-4,3.4)(-3.6,4)(-2.5,4.5)(-1.9,4.4)(-1.4,4)                % top connecting curve
\psecurve(-4,-3.4)(-3.6,-4)(-2.5,-4.5)(-1.9,-4.4)(-1.4,-4)           % borrom connecting curve 

\rput(-1.2,0){$bD\cap U_1$}
\psline[linewidth=0.2pt]{->}(-0.1,0)(0.6,0)

\rput(4.6,0){$bD'_1$}
\psline[linewidth=0.2pt]{->}(4,0)(2.45,0)

\rput(3.7,-2.8){$bD_1$}
\psline[linewidth=0.2pt]{->}(3.1,-2.6)(1.25,-1)

\rput(3.2,3){$U_1$}
\psline[linewidth=0.2pt]{->}(2.8,2.8)(2.3,2.3)

\rput(0.7,5){$V_1$}

\rput(-7.8,0){$D\cap D'_1$}
\psline[linewidth=0.2pt]{->}(-6.6,0)(-4.8,0)

\end{pspicture}
\caption{The domains $D'_1$ and $D_1$}
\label{Fig1}
\end{figure}

Proposition \ref{sprays-exist} furnishes a dominating $h$-spray
$F\colon \bar D\times P\to X$ of class $\cA^r(D)$ with the core section
$F_0=f$. By shrinking its parameter set $P\subset\C^N$ if necessary 
we can insure that properties (ii) and (iii) above are satisfied 
if we replace $f=F_0$ by any section $F_t=F(\cdotp,t)$, $t\in P$. 

By using the special coordinate chart $(W_1,V_1,\Phi_1)$ 
we find an open set $\Omega\subset V_1$
containing $\bar D\cap V_1$ and a holomorphic $h$-spray 
$G\colon \Omega\times P\to X|_{\Omega}=h^{-1}(\Omega)$, with range contained 
in $W_1$, such that the restriction of $G$ to $(\bar D\cap V_1)\times P$
approximates $F$ as close as desired in the $\cC^r$ topology.
(The size of  $\Omega$ will depend on $G$.) 

We wish to use Proposition \ref{gluing-sprays} to glue $F$ and $G$ into
a single $h$-spray. We cannot do this directly since their domains
do not form a Cartan pair. (The domain of  $G$ need not 
contain all of the set $\bar D'_1$ which forms a Cartan pair with $D$.) 
We proceed as follows. If $G$ is sufficiently close to  $F$ in the 
$\cC^r$ topology on the intersection of their domains, 
Lemma 4.4 in \cite{BDF} furnishes a transition map $\gamma(z,t)=z+c(z,t)$ 
of class $\cA^r$ between $F$ and $G$, defined for 
$z \in \bar D\cap \bar D'_1$ and $t$ in a smaller
parameter set $P'\Subset P\subset \C^N$, such that $\gamma$ is 
$\cC^r$ close to $\gamma_0(z,t)=t$ (depending on the closeness of $G$
to $F$), and $F(z,t)=G(z,\gamma(z,t))$ for $z \in \bar D\cap \bar D'_1$ and 
$t\in P'$.

Theorem 3.2 in \cite{BDF}, applied to $\gamma$ on the Cartan pair $(D,D'_1)$, 
furnishes a smaller parameter set $P_1\subset P'$ and maps 
$\alpha(z,t)=t+a(z,t)$ on $\bar D\times P_1$, $\beta(z,t)=t+b(z,t)$
on $\bar D'_1\times P_1$, of class $\cA^r$ and close to $\gamma_0$ 
on their respective domains, such that 
$\gamma(z,\alpha(z,t))=\beta(z,t)$ for $z\in \bar D\cap \bar D'_1$
and $t$ close to $0$.
The spray $F(z,\alpha(z,t))$ is then defined on $\bar D\times P_1$ and is close
to $F$, the spray $G(z,\beta(z,t))$ is defined on $(\bar D'_1\cap \Omega)\cap P_1$, 
and by the construction these sprays agree for $z\in\bar D\cap\bar D'_1 \subset \Omega$.
The central section $f_1$ of the new amalgamated spray, obtained by setting $t=0$,
is then of class $\cA^r$ on $\bar D\cup (\bar D'_1\cap \Omega)$. 

It remains to restrict $f_1$ to a suitably chosen strongly pseudoconvex
domain $D_1\Subset S$ contained in $D\cup (D'_1\cap \Omega)$ and satisfying
the other required properties. We choose $D_1$ such that it agrees with 
$D$ outside of $V_1$, while 
\[
		D\cap V_1 =\phi_1^{-1}(\{z\in B\colon \tau(z) < \e \chi(z) \})
\]
for a small $\e>0$ (Fig.\ \ref{Fig1}). By choosing $\e>0$ sufficiently small 
(depending on $f_1$) we insure that the hypotheses (i)--(iii) are satisfied by
the pair $(D_1,f_1)$. This completes the first step. 

Applying the same procedure to $(D_1,f_1)$ and the second chart $(W_2,V_2,\Phi_2)$
we get the next pair $(D_2,f_2)$.  After $j_0$ steps we find a section 
$f_{j_0}$ over a neighborhood $D_{j_0}$ of $\bar D$ in $S$ 
which approximates $f$ as close as desired in the $\cC^r$ topology on $\bar D$.
Since all steps were made by homotopies, $f_{j_0}|_{\bar D}$ 
is homotopic to $f=f_0$ in $\Gamma_{\cA}^r(D,X)$.
\end{proof}

{\em Proof of Theorem \ref{Riemann-surfaces}.}
Assume that $D$ is a relatively compact connected domain with $\cC^2$ boundary 
in an open Riemann surface $S$ and $f_0\colon \bar D\to Y$ is a map
of class $\cA(D,Y)$ to an $n$-dimensional complex space $Y$.
We wish to prove that $f_0$ can be uniformly approximated by
maps which are holomorphic in open neighborhoods of $\bar D$ in $S$.
(The analogous result for maps of class $\cA^r(D,Y)$ can be 
proved in a similar way.)

We proceed by induction on the dimension of $Y$. Assume that the result 
already holds for complex spaces of dimension $<n=\dim Y$ (this is trivially 
satisfied when $n=1$). If $f_0(\bar D)\subset Y_{sing}$, 
the inductive hypothesis gives approximation of $f_0$ 
by holomorphic maps to $Y_{sing}$.
(It may happen that the image of any nearby holomorphic map
is contained in $Y_{sing}$; see the example in \cite{FW}.)

Suppose now that $f_0(\bar D)\not\subset Y_{sing}$.
The set $\sigma=\{z\in\bar D\colon f_0(z)\in Y_{sing}\}$ is locally 
the common zero set of finitely many functions of class $\cA(D)$;
hence $\sigma$ is a closed subset of $\bar D$ such that $\sigma\cap D$ 
is discrete in $D$ and $\sigma\cap bD$ has linear measure zero in $bD$. 

We shall need the following result analogous to Corollary \ref{sprays-exist2}.

\begin{lemma}
\label{sprays-exist3}
(Notation as above.) 
There exists a spray of maps $f\colon\bar D \times P\to Y$
of class $\cA(D)$ which is dominating on $\bar D\bs \sigma$
and whose core map is  $f_0$.
\end{lemma}

\begin{proof}
As in the proof of Proposition \ref{sprays-exist}
we choose a finite sequence of domains 
$D_0\subset D_1\subset\cdots\subset D_m=D$  with $\cC^2$
boundaries such that $\bar D_0\subset D$, and
for every $j=0,1,\ldots,m-1$ we have $D_{j+1}=D_j\cup B_j$ 
where $B_j$ is a convex bump on $D_j$ (Def.\ \ref{Cartan-pair}).
In addition, the properties of the set $\sigma=f_0^{-1}(Y_{sing})$
imply that the sets $B_j$ and $D_j$ can be chosen such that
the following hold for each $j=0,\ldots, m-1$:
\begin{itemize}
\item[(i)]  $\sigma \cap bD_j \cap D=\emptyset$,   
\item[(ii)] $\sigma\cap bD_j \cap bD$ is contained in the 
relative interior of $bD_j\cap bD$ (in $bD$),
\item[(iii)] $\sigma\cap \bar D_j\cap \bar B_j=\emptyset$.
\end{itemize}

Lemma 4.2 in \cite{BDF} furnishes sprays $f \colon \bar D_0\times P\to Y$
and $f' \colon \bar B_0\times P \to Y$, satisfying the required properties 
over $\bar D_0$ resp.\ $\bar B_0$. In particular, the core of each of these 
two sprays is $f_0$, restricted to the respective domain, and the
exceptional set (the set where the spray fails to dominate) is $\sigma$.
Note that both sprays are dominating over the convex set 
$\bar C_0=\bar D_0\cap \bar B_0$ by (iii), and they agree at $t=0$.

The proof of Lemma 4.4 in \cite{BDF}  gives a map $\gamma(z,t)$ 
of class $\cA(C_0 \times P)$ (after shrinking $P$ around $0$),
satisfying $f(z,t)=f'(z,\gamma(z,t))$ for $(z,t)\in \bar C_0\times P$
and $\gamma(z,0)=0$ for $z\in \bar C_0$. (In that proof we have assumed
that the two sprays are close to each other, but in the
present situation this is not necessary since $f(z,0)=f'(z,0)$ 
for $z\in \bar C_0$. It suffices to find a direct summand of class
$\cA(C_0)$ in $\bar D_0\times\C^N$ to the kernel of the $t$-derivative 
of each of the two sprays and apply the implicit function theorem
as in \cite{BDF}.) 

As in the proof of Lemma \ref{approximation1} we approximate
$\gamma$ uniformly on $\bar C_0\times P$ by an entire map
$\gamma'(z,t)$ satisfying $\gamma'(z,0)=0$ for $z\in\C^n$. 
The spray $(z,t)\to f'(z,\gamma'(z,t))$ is then 
defined on $(\bar C_0\times P)\cup (\bar B_0 \times P_0)$
for some polydisc $P_0\subset\C^n$ around $0$ 
(which may depend on $\gamma'$), and it agrees with $f$ for $t=0$. 
If the two sprays are sufficiently uniformly close 
to each other on $\bar C_0\times P$ (as we may assume to be the case), 
we can glue them into a new spray $\wt f$ with the required properties 
over $\bar D_1=\bar D_0\cup \bar B_0$.  
(See the proof of Proposition \ref{sprays-exist} 
above, or the proof of Proposition 4.3 in \cite{BDF}. 
The important point is that our sprays agree along $t=0$,
so we get a nonempty parameter set for the new spray.)  
Note also that the gluing process does not increase the 
exceptional set of the spray. After finitely many steps of 
this kind we get a desired spray $f$ over $\bar D$.
\end{proof}

To complete the proof of Theorem \ref{Riemann-surfaces}
we proceed exactly as in the proof of Theorem 5.1 in \cite{BDF}.
The main idea is to attach to $\bar D$ a convex bump $B$
such that $\sigma\cap \bar D\cap \bar B =\emptyset$
(this is possible since $\sigma\cap bD$ has empty interior in $bD$)
and then approximate $f$ by another spray over 
$K= \bar D\cup \bar B$. Let $g\colon K \to Y$ denote 
the core of the new spray. There exists a holomorphic vector field
$\xi$ in a neighborhood of $K$ in $S$ which points into $D$
at every point of $bD\bs B$ (we put no condition on $\xi$
at the points of $bD\cap \bar B$; see \cite{BDF}). 
Denoting its flow by $\phi_t$, the map $g\circ \phi_t$ 
is defined and holomorphic in an open neighborhood of $\bar D$
for every sufficiently small $t>0$, and it approximates
$g$ (and hence $f_0$) uniformly on $\bar D$.
\qed
\smallskip

Combining Theorem \ref{approximation-sections} with the main result of \cite{ANN}
gives the following version of the Oka principle.
Note that Corollary \ref{cor1} is just Corollary \ref{improved-Oka} 
applied to the product submersion $X=S\times Y\to S$.

\begin{corollary}
\label{improved-Oka}
Assume that $D\Subset S$ is a strongly pseudoconvex domain 
with $\cO(S)$-convex closure in a Stein manifold $S$
and $h\colon X\to S$ is a holomorphic fiber bundle whose fiber $Y$
enjoys {\rm CAP}  (Def.\ \ref{CAP}). For any continuous section $f_0\colon S\to X$ 
which is holomorphic on $D$ there exists a homotopy $f_t\colon S\to X$ 
of continuous sections such that $f_t|_D$ is holomorphic and uniformly 
close to $f_0|_D$ for each $t\in [0,1]$, and $f_1$ is holomorphic on $S$. 

Furthermore, every homotopy of continuous sections $f_t\colon S\to X$, with
$f_t|_D$ holomorphic for each $t\in [0,1]$ and $f_0,f_1$ holomorphic on $S$, 
can be deformed with fixed ends to a homotopy $\wt f_t\colon S\to X$ 
consisting of sections which are holomorphic on $S$ such that the 
entire homotopy remains holomorphic on~$D$. 

The analogous results hold for sections of class $\cC^r$ provided
that $bD\in \cC^\ell$ $(\ell\ge 2)$ and $r\in \{0,1,\ldots,\ell\}$.
\end{corollary}

\begin{proof}
By Theorem \ref{approximation-sections} we can approximate $f$ as close as 
desired in the $\cC^r(\bar D,X)$ topology by a $\cC^r$ section $g\colon U\to X$
which is holomorphic in an open neighborhood $U$ of $\bar D$ in $S$.
If the approximation of $f$ by $g$ is sufficiently close and $U$ is chosen sufficiently 
small then $g$ is homotopic to $f|_U$ by a homotopy 
of $\cC^r$ sections $g_t\colon U\to X$ $(t\in[0,1],\ g_0=f,\ g_1=g)$ 
which are holomorphic in $D$ (Corollary \ref{TNT}). 
Choose a smooth function $\chi\colon S\to [0,1]$ with compact support 
contained in $U$ such that $\chi=1$ in a smaller open neighborhood of $\bar D$.
Set $\tilde g_t(z) = g_{\chi(z) t}(z)$; this is a homotopy of $\cC^r$ sections
which are holomorphic on $D$, they extend to all of $S$ and equal $f$ on $S\bs U$.

Assuming that the fiber $Y$ of $h\colon X\to S$ enjoys CAP,
Theorem 1.2 in \cite{ANN} shows that $\tilde g_1$ is homotopic 
to a holomorphic section $f_1\colon S\to X$ by a homotopy 
of sections $h_t\colon S\to X$ $(t\in[0,1],\ h_0= \tilde g_1,\ h_1=f_1)$
which are holomorphic and uniformly close to $\tilde g_1$ 
in an open neighborhood of $\bar D$. By combining the homotopies 
$\tilde g_t$ and $h_t$ we obtain a homotopy from 
$f=f_0$ to $f_1$ satisfying the stated properties.

Similarly, a homotopy $\{f_t\}_{t\in[0,1]}$ in the second statement can be deformed 
(with fixed ends at $t=0,1$) to a homotopy 
$\{f'_t\}_{t\in[0,1]}$ consisting of sections which are holomorphic 
in an open neighborhood of $\bar D$ in $S$;  it remains to apply 
the one-parametric Oka principle \cite[Theorem 5.1]{ANN} to $\{f'_t\}$.
\end{proof}

\section{The Oka property for $\cA^r(D)$-bundles with flexible fibers}
Let $D$ be a relatively compact strongly pseudoconvex domain in 
a Stein manifold $S$. Recall that $\cA(D)$ denotes the algebra of all continuous 
functions on $\bar D$ which are holomorphic on $D$;
we shall use the analogous notation for (open or closed) domains  in $\bar D$.

A compact set $K$ in $\bar D$ is {\em $\cA(D)$-convex} if for every $p\in \bar D\bs K$ 
there is a $g\in\cA(D)$ satisfying $|g(p)|>\sup_{z\in K} |g(z)|$.

The following result is a precise version of Theorem \ref{main}.

%
%
%  Theorem main-bis
%
%
\begin{theorem} 
\label{main-bis}
Assume that $S$ is a Stein manifold, $D\Subset S$ is a strongly pseudoconvex domain 
with $\cC^\ell$ boundary $(\ell\ge 2)$, $r\in\{0,1,\ldots,\ell\}$, and
$h\colon X\to \bar D$ is an $\cA^r_Y(D)$-bundle (Def.\ \ref{Ar-bundle}).
Choose a distance function $d$ on the manifold $J^r(\bar D,X)$ of all $r$-jets 
of sections $\bar D\to X$ of $h$. Let $K$ be a compact $\cA(D)$-convex subset of $\bar D$ 
and let $U\subset S$ be an open set containing $K$. 
If the fiber $Y$ enjoys {\rm CAP} (Def.\ \ref{CAP})
then sections $\bar D\to X$ satisfy the following:
\begin{itemize}
\item[(i)] Given a continuous section $f_0\colon \bar D\to X$ which is of class $\cC^r$ 
on $U\cap \bar D$ and holomorphic in $U\cap  D$, there exist for every $\e>0$ 
an open neighborhood $V\subset U$ of $K$ and a homotopy of sections 
$f_t\colon \bar D\to X$ $(t\in[0,1])$ which are of class $\cC^r$ on $V\cap \bar D$ 
and holomorphic on $V\cap D$ such that $f_1\in \Gamma_\cA^r(\bar D,X)$ and
\[
		\sup_{x\in K} d\bigl(j^r_x f_t,j^r_x f_0\bigr) <\e,\quad t\in[0,1].
\]
\item[(ii)] Given a homotopy of continuous sections $f_t\colon \bar D\to X$ $(t\in[0,1])$
such that $f_0,f_1\in\Gamma_\cA^r(\bar D,X)$, 
$f_t$ is holomorphic in $U\cap  D$ and of class $\cC^r$ on $U\cap \bar D$ 
for each $t\in[0,1]$, with continuous dependence of $j^r f_t|_{U\cap \bar D}$ on $t$,
there are an open neighborhood $V\subset U$ of $K$ and a homotopy of sections 
$g_{t,s} \colon \bar D\to X$ $(t,s\in[0,1])$ satisfying 
\smallskip
\begin{itemize}
\item[(1)] $g_{t,0}=f_t$, $g_{0,s}=f_0$, $g_{1,s}=f_1$ for all  $t,s\in[0,1]$,
\item[(2)] $g_{t,1}\in \Gamma_\cA^r(\bar D,X)$ for all $t\in [0,1]$,
\item[(3)] $g_{t,s}$ is of class $\cC^r$ on $V\cap \bar D$, holomorphic on 
$V\cap D$, and 
\[
		\sup_{x\in K} d\bigl(j^r_x g_{t,s},j^r_x f_t\bigr) <\e, \quad s,t\in[0,1].
\]
\end{itemize}
\end{itemize}
\end{theorem}

\begin{proof}
The proof of (i) will consist of two parts. In the first part we 
find a section $g_1$ of $X$ satisfying (i) 
on a strongly pseudoconvex domain $\bar D_1\subset \bar D$
containing $K$ and such that $\bar D$ is obtained from $\bar D_1$ 
by attaching finitely many convex bumps. In part 2 we show how to 
approximately extend the solution over each bump 
to get a solution $f_1$ on all of $\bar D$.

{\em Part 1.}  We have $D=\{\rho <0\}$ where $\rho$
is a $\cC^\ell$ strongly plurisubharmonic function in a neighborhood
of $\bar D$, with $d\rho\ne 0$ on $bD=\rho^{-1}(0)$.

Consider first the case when $K\subset D$. Choose
a constant $c_1<0$ close to $0$ such that 
$K\subset D_1:= \{z\in D\colon \rho(x) <c_1\}$ 
and $d\rho\ne 0$ on $\{c_1\le \rho \le 0\}= \bar D\bs D_1$.
Since the fiber $Y$ enjoys CAP, the Oka principle on open Stein
manifolds \cite[Theorem 1.2]{ANN} gives a holomorphic section 
$g_1$ of $X$ over an open neighborhood of $\bar D_1$ 
which satisfies the conclusion (i) over $\bar D_1$
(with $f_1$ replaced by $g_1$). Now proceed directly to Part 2 below. 

The situation is more complicated when $K\cap bD\ne \emptyset$.
Let $D_0=\{\rho<c_0\}$ for some $c_0<0$ close to $0$ such that 
$d\rho\ne 0$ on $\{c_0\le \rho \le 0\}$.  If $c_0$ is chosen sufficiently 
close to $0$, there is a strongly pseudoconvex domain $B_0\subset D$ 
such that $K\subset \bar B_0$, $(B_0,D_0)$ is a Cartan pair of class $\cC^\ell$,
and $\bar B_0\cap \bar D_0$ is holomorphically convex in $D$.
(Here we used the assumption that $K$, and hence $K\cap \bar D_0$,
is $\cA(D)$-convex.) Set $D_1=D_0\cup B_0$ (Fig.\ \ref{FigX}). 
Suitable choices of $D_0$ and $B_0$ insure that 
$\bar D$ is obtained from $\bar D_1$ by attaching finitely 
many small convex bumps \cite[Lemma 12.3]{HL0}.

\begin{figure}[ht]
\psset{unit=0.4 cm} 
\begin{pspicture}(-9,-7)(9,7)

%\psgrid
%
%  large circle
%
\pscircle[linecolor=DarkBlue](0,0){6.06}       
\rput(-8,0){$bD$}
\psline[linewidth=0.2pt]{->}(-7.2,-0.1)(-6.2,-0.1)

%
%  small circle
%
\pscircle[linecolor=red](0,0){5.3}       
\rput(-7.3,-3){$bD_0$}
\psline[linewidth=0.2pt]{->}(-6.5,-2.8)(-4.95,-1.95)

%
%  domain B_0
%
\pscustom[fillstyle=hlines,hatchcolor=yellow]{
\psccurve[linecolor=Black](2.6,4.6)(5.4,2.6)(6,0)(5.4,-2.6)(2.6,-4.6) (1,-5)(-1,0)(1,5)(2.6,4.6)  
%\pscurve[linecolor=Purple](2.6,-4.6)(1,-5)(-1,0)(1,5)(2.6,4.6) 
}
\rput(1.4,-4){$B_0$}

%
%
%  The set U
%       
\pscurve[linecolor=Black,linestyle=dashed](1.2,5.9)(-1,3)(-1.7,0)(-1,-3)(1.2,-5.9)
%\rput(4.2,6){$bU\cap D$}
%\psline[linewidth=0.3pt]{->}(2.6,5.9)(0.9,5.6)

\psline[linewidth=0.2pt]{<-}(-1.1,2.6)(0.5,2.6)
\rput(1.2,2.7){$bU$}

%
%
% K
%
\psarc(0,0){6.06}{340}{20}
\psecurve[linecolor=Black](5.8,1.8)(5.7,2)(4,2.6)(3,2.2)(2,1.4)(0.4,0.4)
(1,-1.8)(2,-3)(3,-2.4)(5.7,-2)(5.8,-1.8)
\rput(2.5,0){$K$}
\end{pspicture}
\caption{The sets $K$, $D_0$, $B_0$ and $U$}
\label{FigX}
\end{figure}
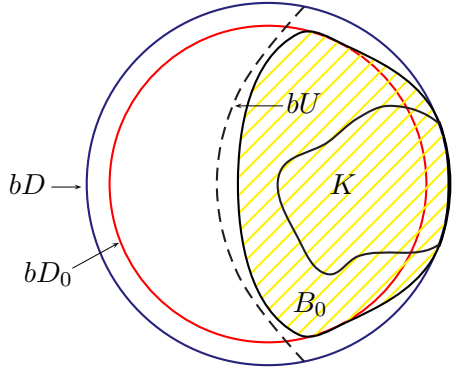

Let $U$ be as in the statement of the theorem.
Choose a strongly pseudoconvex domain $D' \subset D\cap U$ 
such that $\bar B_0 \subset D'\cup bD$. By Corollary \ref{sprays-exist2} 
there exist a domain $P_0\subset \C^N$ containing the origin and a 
dominating $h$-spray $f\colon \bar D'\times P_0\to X$ 
of class $\cA^r(D')$ with the central section $f_0|_{\bar D'}$.

Let $P\Subset P_0$ be a ball centered at the origin in $\C^N$.
Since the set $(\bar B_0\cap \bar D_0)\times \bar P$ is holomorphically
convex in $D\times P_0$ and $Y$ enjoys CAP, the Oka principle
for open Stein manifolds \cite[Theorem 1.2]{ANN} gives an
open neighborhood $W\subset D$ of $\bar D_0$ and a holomorphic 
spray of sections $f'\colon W\times P_0 \to X$ which approximates
$f$ uniformly on a neighborhood of $(\bar B_0\cap \bar D_0)\times \bar P$;
by the Cauchy estimates the approximation is then $\cC^r$ close 
on the latter set. Now we apply Proposition \ref{gluing-sprays}
to glue the spray $f$ (over $\bar B_0$) and the spray $f'$ 
(over $\bar D_0$) to a new spray $g\colon \bar D_1\times P' \to X$
of class $\cA^r(D_1)$ which approximates $f$ in the $\cC^r$ topology
on $(\bar B_0\cap \bar D_0)\times P'$. Its central section 
$g_1=g(\cdotp,0)\colon \bar D_1\to X$ then satisfies condition 
(i) over $\bar D_1$.

{\em Part 2:}
By the choice of $D_1$ there exist strongly pseudoconvex domains 
$D_1\subset D_2\subset \cdots\subset D_m = D$ with $\cC^\ell$ boundaries 
such that for every $j=1,\ldots,m-1$ we have $D_{j+1}=D_j\cup B_j$, 
where $B_j$ is a convex bump on $D_j$ (Def.\ \ref{Cartan-pair}).
Each of the sets $B_1,\ldots,B_{m-1}$ is chosen so small that 
the fiber bundle $X$ is trivial over its closure.
To find a solution $f_1\colon \bar D \to X$ we successively extend 
the section $g_1\colon \bar D_1\to X$ found in Part 1 over each 
bump by using a finite induction; this will produce 
sections $g_j\in \Gamma^r_\cA(D_j,X)$ $(j=2,3,\ldots,m)$
such that the restriction of $g_j$ to $\bar D_{j-1}$ is homotopic to
$g_{j-1}$, and is close to $g_{j-1}$ in $\Gamma_\cA^r(\bar D_{j-1},X)$.
Since all steps are of the same kind, it suffices to explain
how to obtain $g_2$ from $g_1$. 

Corollary \ref{sprays-exist2} furnishes a ball $P \subset\C^N$
containing the origin and a dominating $h$-spray 
$G_1\colon \bar D_1\times P \to X$ 
of class $\cA^r(D_1)$ with the central section $g_1$. 
Let $P'\Subset P$ be a smaller ball around the origin.
Since the fiber bundle $X$ is trivial over $\bar B_1$, we can 
identify sections of $X$ over (subsets of) $\bar B_1$ with maps 
to the fiber $Y$. Also, using local holomorphic coordinates in 
a neighborhood of $\bar B_1$ in $S$ we shall identify 
$\bar B_1\cap \bar D_1 \subset \bar B_1$ with 
compact convex sets in $\C^n$. 

Using these identifications, we first approximate $G_1$
in the $\cC^r$ topology on $(\bar B_1\cap \bar D_1) \times \bar P'$
by a holomorphic map from an open neighborhood of this set (in $\C^n\times\C^N$)
to $Y$ (just precompose $G_1$ with linear contractions to an interior point
of its domain);  by using the CAP property of $Y$ we then approximate this 
map by an entire map $\C^n\times\C^N \to Y$. Passing back to subsets of $S$, 
the above procedure gives a spray 
$G'_1 \colon \bar B_1\times P' \to X$ of class $\cA^r(B_1)$ 
which approximates $G_1$ as close as desired in the 
$\cC^r$ topology on $(\bar D_1\cap \bar B_1)\times P'$.
Finally we glue $G_1$ and $G'_1$ into a spray 
$G_2\colon \bar D_2\times P'' \to X$ of class $\cA^r(D_2)$ 
by appealing to Proposition \ref{gluing-sprays}.
Letting $g_2=G_2(\cdotp,0) \colon \bar D_2\to X$ 
be the central section of $G_2$ we complete the induction step.

After $m$ steps of this kind we obtain a section 
$g_m \in \Gamma_\cA^r(D,X)$ (which we now call $f_1$) 
satisfying the conclusion of Theorem \ref{main-bis} (i). 
A homotopy from $f_0$ to $f_1$ with the stated properties is obtained by 
combining the homotopies obtained in the individual steps of the proof.

Part (ii) is proved in essentially the same 
way by using the equivalence between CAP and the one-parametric Oka property
\cite[\S 5]{ANN}. A one-parametric spray with the given center
is furnished by the parametric part of Corollary \ref{sprays-exist2}, 
and we glue sprays by the parametric version of Proposition \ref{gluing-sprays}.
\end{proof}

In conclusion we mention the following (apparently) open problem.

\begin{problem} 
Assume that $D$ is a smoothly bounded strongly pseudoconvex 
domain in a Stein manifold, $r\in \N\cup\{\infty\}$, and $X\to\bar D$ is a fiber 
bundle of class $\cA^r(D)$. Let $f\colon \bar D\to X$ be a section of class $\cA(D)$.
Is it possible to approximate $f$ uniformly on $\bar D$ by sections
of class $\cA^r(D)$~?
\end{problem}

\textit{Added in the final revision.} 
In the subsequent paper \cite{FF:manifolds} of the second author
the conclusion of Theorem \ref{Banach-manifold} is extended 
to H\"older and Sobolev spaces of holomorphic maps, using the fact 
that the graph 
of any map $f\colon\bar D\to Y$ of class $\cA(D)$ admits a basis of 
open Stein neighborhoods in the ambient manifold $S\times Y$
\cite[Theorem 1.2]{FF:manifolds}. However, open Stein neighborhoods
cannot be used in the proofs of Theorems \ref{main} and \ref{main-bis} 
since the bundle is only defined over $\bar D$.

%
%
%  THANKS, THANKS AND THANKS
%
%
%

\smallskip
\textit{Acknowledgements.} 
We wish to thank J.\ Leiterer for pointing out the works of 
Heunemann and Sebbar, and L.\ Lempert for helpful remarks
concerning the theory of complex Banach manifolds.

\bibliographystyle{amsplain}

\end{document}